\documentstyle[12pt]{article}
\begin{document}
\newtheorem{theorem}      {Th\'eor\`eme}[section]
\newtheorem{theorem*}     {theorem}
\newtheorem{proposition}  [theorem]{Proposition}
\newtheorem{definition}   [theorem]{Definition}
\newtheorem{e-lemme}        [theorem]{Lemma}
\newtheorem{cor}   [theorem]{Corollaire}
\newtheorem{resultat}     [theorem]{R\'esultat}
\newtheorem{eexercice}    [theorem]{Exercice}
\newtheorem{rrem}    [theorem]{Remarque}
\newtheorem{pprobleme}    [theorem]{Probl\`eme}
\newtheorem{eexemple}     [theorem]{Exemple}
\newcommand{\preuve}      {\paragraph{Preuve}}
\newenvironment{probleme} {\begin{pprobleme}\rm}{\end{pprobleme}}
\newenvironment{remarque} {\begin{rremarque}\rm}{\end{rremarque}}
\newenvironment{exercice} {\begin{eexercice}\rm}{\end{eexercice}}
\newenvironment{exemple}  {\begin{eexemple}\rm}{\end{eexemple}}
%
%
\newtheorem{e-theo}      [theorem]{Theorem}
\newtheorem{theo*}     [theorem]{Theorem}
\newtheorem{e-pro}  [theorem]{Proposition}
\newtheorem{e-def}   [theorem]{Definition}
\newtheorem{e-lem}        [theorem]{Lemma}
\newtheorem{e-cor}   [theorem]{Corollary}
\newtheorem{e-resultat}     [theorem]{Result}
\newtheorem{ex}    [theorem]{Exercise}
\newtheorem{e-rem}    [theorem]{Remark}
\newtheorem{prob}    [theorem]{Problem}
\newtheorem{example}     [theorem]{Example}
\newcommand{\proof}         {\paragraph{Proof~: }}
\newcommand{\hint}          {\paragraph{Hint}}
\newcommand{\heuristicproof}{\paragraph{heuristic proof}}
\newenvironment{e-probleme} {\begin{e-pprobleme}\rm}{\end{e-pprobleme}}
\newenvironment{e-remarque} {\begin{e-rremarque}\rm}{\end{e-rremarque}}
\newenvironment{e-exercice} {\begin{e-eexercice}\rm}{\end{e-eexercice}}
\newenvironment{e-exemple}  {\begin{e-eexemple}\rm}{\end{e-eexemple}}
\newcommand{\1}        {{\bf 1}}
\newcommand{\pp}       {{{\rm I\!\!\! P}}}
\newcommand{\qq}       {{{\rm I\!\!\! Q}}}
\newcommand{\B}        {{{\rm I\! B}}}
\newcommand{\cc}       {{{\rm I\!\!\! C}}}
\newcommand{\N}        {{{\rm I\! N}}}
\newcommand{\R}        {{{\rm I\! R}}}
\newcommand{\D}        {{{\rm I\! D}}}
\newcommand{\Z}       {{{\rm Z\!\! Z}}}
\newcommand{\C}        {{\bf C}}        
\newcommand{\rank}{\hbox{rank}}
\newcommand{\CC}{{\cal C}}
\def \Re {{\rm Re\,}}
\def \Im {{ \rm Im\,}}
\def\Hom{{\rm Hom\,}}
\def\Lip{{\rm Lip}}
\def\st{{\rm st}}
\def\J{{\cal J}}
\def\A{{\cal A}}
\def\Mat{{\rm Mat}}
\def\m{{\rm m}}
\def\ind{{\rm ind\,}}
\def\<{\langle}
\def\>{\rangle}
\def\bar{\overline}
%
%
\newcommand{\dontforget}[1]
{{\mbox{}\\\noindent\rule{1cm}{2mm}\hfill don't forget : #1
\hfill\rule{1cm}{2mm}}\typeout{---------- don't forget : #1 ------------}}
\newcommand{\note}[2]
{ \noindent{\sf #1 \hfill \today}

\noindent\mbox{}\hrulefill\mbox{}
\begin{quote}\begin{quote}\sf #2\end{quote}\end{quote}
\noindent\mbox{}\hrulefill\mbox{}
\vspace{1cm}
}
\title{Regularization of almost complex structures
and gluing holomorphic discs to tori}
\author{ Alexandre Sukhov{*} and Alexander Tumanov{**}}
\date{}
\maketitle

{\small

* Universit\'e des Sciences et Technologies de Lille, Laboratoire
Paul Painlev\'e,
U.F.R. de
Math\'e-matique, 59655 Villeneuve d'Ascq, Cedex, France,
 sukhov@math.univ-lille1.fr

** University of Illinois, Department of Mathematics
1409 West Green Street, Urbana, IL 61801, USA, tumanov@illinois.edu
}
\bigskip

Abstract. We prove a result on removing singularities of almost
complex structures pulled back by a non-diffeomorphic map.
As an application we prove the existence of
global $J$-holomorphic discs with boundaries attached
to real tori.
\bigskip


\bigskip

\section{Introduction}
In this paper we prove a general result (Theorem \ref{removal})
on removing singularities of almost complex structures pulled
back by a non-diffeomorphic map.
In our joint paper  with Bernard Coupet \cite{CoSuTu}
we use special coordinates in an almost complex manifold $(M,J)$
to reduce a boundary value problem for $J$-holomorphic
discs in $M$ to that for quasi-linear PDE in the plane.
In \cite{CoSuTu} we make a simplifying assumption that the
coordinates are introduced by a locally diffeomorphic map,
although it is not the case in general.
The main result of the present work allows to
extend the range of applications of the methods of \cite{CoSuTu}
to the non-diffeomorphic case.
We use Theorem \ref{removal} in order to obtain  results on
attaching  pseudo-holomorphic discs to real tori.
These results are new even in the complex Euclidean
space $\cc^2$. Thus the use of almost complex structures
leads to new results in the classical complex analysis.

We now describe the main results and the organization of the paper.
In Section 2 we recall the notion of an almost complex structure $J$,
in particular, we include some details on representing $J$
by a complex matrix function.
In Section 3 we prove Theorem \ref{removal} mentioned above.
It deals with a map from the standard bidisc with coordinates
$(z,w)$ to an almost complex manifold $(M,J)$.
The map takes the coordinate lines $z=c$ to a given family
of $J$-holomorphic discs.
The map is not necessarily locally diffeomorphic.
Nevertheless, we prove that under natural assumptions,
the pullback of $J$ exists and is sufficiently regular:
H\"older in $z$ and Lipschitz in $w$.
The result is useful even for integrable $J$;
in this case it also holds in higher dimension.
As a side production of the proof, we obtain
results (Propositions \ref{Holder} and \ref{final})
on regularity of generalized analytic functions
in the sense of Vekua \cite{Ve} and a result
(Theorem \ref{interpolation1}) on decomposition
of the phase function of a complex polynomial.
We hope that Theorem \ref{removal} will find further
applications, in particular, in the theory of foliations.

Theorem \ref{removal} is relevant to blow-up situations.
Theory of blow-ups for almost complex manifolds is not yet
fully developed. We would like to mention a result
by Duval \cite{Du1} in complex dimension 2, in which
after a blow-up the resulting structure is no longer smooth.
We hope that our results will be useful in developing
blow-up techniques in almost complex category.

In Section 4 we construct pseudo-holomorphic discs
attached to the standard torus in $\cc^2$ equipped
with a certain almost complex structure. This improves
the corresponding result of \cite{CoSuTu} by adding
a continuous depending statement to it.
In Section 5 we prove existence theorems for pseudo-holomorphic
discs (Theorems \ref{tori} and \ref{toridouble}) with given boundary
conditions. Under various assumptions, we construct
pseudo-holomorphic discs with boundary glued to real tori,
which are not Lagrangian in general.

In his pioneering work, Gromov \cite{Gr} proved
the existence of pseudo-holomorphic discs glued
to smooth Lagrangian submanifolds.
Ivashkovich and Shevchishin \cite{IvSh1}
extended the result to the case of immersed
Lagrangian submanifolds.
Forstneri\v c \cite{Fo} constructed discs attached
to certain totally real 2-tori in the space $\cc^2$ with
the standard complex structure.
Cerne \cite{Ce} generalized the result of \cite{Fo}
to the case of bordered Riemann surfaces.
On the other hand,
Alexander \cite{Al} constructed a
totally real 2-dimensional torus in $\cc^2$ which contains no
boundary of a holomorphic disc.
Moreover, Duval \cite{Du} gave  an example of a torus
with the same property and which in addition is contained
in the unit sphere in $\cc^2$.
Thus some restrictions on the geometry of a torus are necessary
for gluing holomorphic discs to it. We stress that no Lagrangian
conditions are required in Theorem \ref{tori} so our approach
provides a new tool for constructing global pseudo-holomorphic
discs with prescribed boundary conditions.

Abundance and flexibility of real changes of coordinates
allowed in Theorem \ref{removal} represent a contrast
with the rigidity of holomorphic maps.
The flexibility comes at a price because
the Cauchy--Riemann equations for pseudo-holomorphic discs
are non-linear.
However, this analytic difficulty can be handled by
the general theory of elliptic PDE in the plane.
We hope that our methods will find other applications.

We thank the referee for many useful remarks letting
us improve the quality of the paper.
In particular, the referee provided us with a simple
proof of Theorem \ref{removal} for integrable structures,
which we include in Appendix.

A large part of the work was done when the first author was visiting
the University of Illinois in the spring of 2008. He thanks
this university for support and hospitality.

\section{Almost complex manifolds}

Let $(M,J)$ be a smooth almost complex manifold.
Denote by $\D$  the
unit disc in $\cc$ and by $J_{\st}$  the standard complex structure
of $\cc^n$; the value of $n$ is usually clear from the context.
Recall that a smooth map $f:\D \to M$
is called {\it $J$-holomorphic}  if $df \circ J_{\st} = J \circ
df$. We also call such a map $f$ a $J$-{\it holomorphic} disc,
a $J$-disc, a {\it pseudo-holomorphic} disc, or
a {\it holomorphic} disc if $J$ is fixed.

An important result due to  Nijenhuis and Woolf \cite{NiWo}
states that for a given point $p \in M$ and a tangent
vector $v \in T_pM$ there exists a $J$-holomorphic disc
$f: \D \to M$ such that $f(0)=p$ and
$df(0)(\frac{\partial}{\partial\xi}) = \lambda v$ for
some $\lambda> 0$. Here $\xi+i\eta=\zeta\in\cc$.
The disc $f$ can be chosen smoothly depending on the
initial data $(p,v)$ and the structure $J$.

In local coordinates $z\in\cc^n$, an almost complex structure
$J$ is represented by a $\R$-linear operator
$J(z):\cc^n\to\cc^n$, $z\in \cc^n$ such that $J(z)^2=-I$,
$I$ being the identity. Then the Cauchy-Riemann equations
for a $J$-holomorphic disc $z:\D\to\cc^n$ have the form
$$
z_\eta=J(z)z_\xi,\quad
\zeta=\xi+i\eta\in\D.
$$
Following Nijenhuis and Woolf \cite{NiWo},
we represent $J$ by a complex $n\times n$
matrix function $A=A(z)$ so that the Cauchy-Riemann equations
have the form
\begin{eqnarray}
\label{holomorphy}
z_{\bar\zeta}=A(z)\bar z_{\bar\zeta},\quad
\zeta\in\D.
\end{eqnarray}
We first discuss the relation between $J$ and $A$ for
fixed $z$. Let $J:\cc^n\to\cc^n$ be a $\R$-linear map
so that $\det(J_\st+J)\ne0$, here $J_\st v=iv$.
Put
\begin{eqnarray}
\label{mapQ}
Q=(J_\st+J)^{-1}(J_\st-J).
\end{eqnarray}
\begin{e-lemme}
\label{antilinear}
$J^2=-I$ if and only if
$QJ_\st+J_\st Q=0$, that is, $Q$ is complex anti-linear.
\end{e-lemme}
\proof
Put $K=J_\st J$.
Then (\ref{mapQ}) is equivalent to
\begin{eqnarray}
\label{mapK}
Q=(I-K)^{-1}(I+K).
\end{eqnarray}
Note that $(I-K)^{-1}$ and $I+K$ commute.
Then $QJ_\st+J_\st Q=0$ is equivalent to
$(I+K)J_\st(I-K)+(I-K)J_\st(I+K)=0$.
Now using $J_\st^2=-I$ and $K=J_\st J$ we obtain $J^2=-I$.
The lemma is proved.
\medskip

We introduce
$$
\J=\{J:\cc^n\to\cc^n: J\;{\rm is}\;\R{\rm-linear},\;J^2=-I,\;
\det(J_\st+J)\ne0  \}
$$
$$
\A=\{A\in\Mat(n,\cc): \det(I-A\bar A)\ne0 \}
$$

Let $J\in\J$. Then by Lemma \ref{antilinear}, the map $Q$
defined by (\ref{mapQ}) is anti-linear, hence, there is a unique
matrix $A \in\Mat(n,\cc)$ such that
\begin{eqnarray}
\label{eqA}
Av=Q\bar v, \quad
v\in\cc^n.
\end{eqnarray}
The following result essentially is contained in \cite{NiWo}.
\begin{e-pro}
\label{JtoA}
The map $J\mapsto A$ is a birational
homeomorphism $\J\to\A$.
\end{e-pro}
\proof
We first note that if $Q$ is anti-linear, then
$1$ is an eigenvalue of $Q$ if and only if $-1$
is an eigenvalue of $Q$, which in turn holds if and only
if $1$ is an eigenvalue of the complex linear map
$Q^2$. In fact, $Q^2=A\bar A$ because by (\ref{eqA}) we have
$A\bar Av=A\bar{A\bar v}=A\bar{Qv}=Q^2v$.
Hence $Q$ has eigenvalues $\pm1$ if and only if
$A\bar A$ has eigenvalue $1$.

Let $J\in\J$. We show that $A\in\A$.
We again use $K=J_\st J$.
Note that $\det(J_{st}+J)=0$ if an only if
$1$ is an eigenvalue of $K$.
We claim that $Q$ defined by (\ref{mapK})
does not have eigenvalue $-1$. Indeed,
$Qv=-v$ implies $(I+K)v=-(I-K)v$ and $v=0$.
Hence $1$ is not an eigenvalue of $A\bar A$,
that is, $A\in\A$.

Conversely, given $A\in\A$, we show that
there exists a unique $J\in\J$, such that
$J\mapsto A$.
Define $Q$ by (\ref{eqA}).
Then $Q^2=A\bar A$ does not have eigenvalue $1$,
hence $Q$ does not have eigenvalue $-1$.
Then we can find $K$ from (\ref{mapK}) which
yields
$K=-(I+Q)^{-1}(I-Q)$.
This implies that $1$ is not an eigenvalue of
$K$ in the same way that (\ref{mapK}) implies
$-1$ is not an eigenvalue of $Q$.
Define $J=-J_{st}K$. Then $\det(J_{st}+J)\ne0$.
Since $Q$ is anti-linear, then by Lemma \ref{antilinear},
we have $J^2=-I$. Thus $J\in\J$.
The proposition is now proved.
\medskip

The above proof yields a useful formula
of $J$ in terms of $A$ that we include for
future references. Since
$(I+Q)(I-Q)=I-Q^2=I-A\bar A$, then
$(I+Q)^{-1}=(I-A\bar A)^{-1}(I-Q)$.
Hence
$J=-J_{st}K=J_{st}(I-A\bar A)^{-1}(I-Q)^2
=J_{st}(I-A\bar A)^{-1}(I+A\bar A-2Q)$.
Finally,
$$
Jv=i(I-A\bar A)^{-1}[(I+A\bar A)v-2A\bar v].
$$

Let $J$ be an almost complex structure in a domain
$\Omega\subset\cc^n$. Suppose $J(z)\in\J$, $z\in\Omega$.
Then by Proposition \ref{JtoA}, $J$ defines
a unique complex matrix function $A$ in $\Omega$ such that
$A(z)\in\A$, $z\in\Omega$. We call $A$ the
{\it complex matrix} of $J$. The matrix $A$ has the
same regularity properties as $J$.

\section{Removing singularities of almost complex structures}

Our construction of discs with prescribed boundary conditions is
based on a suitable choice of coordinate systems. As we
will see later, it is useful for applications to allow
changes of coordinates which are not necessarily locally
diffeomorphic. This presumably leads to singularities of
almost complex structures obtained by non-diffeomorphic
changes of coordinates. Under some mild assumptions, we prove
that such singularities are removable.

As usual, we denote by $C^{k,\alpha} \, (k\ge 0,\; 0<\alpha\le1)$
the class of functions whose derivatives to order $k$
satisfy a H\"older (Lipschitz) condition with exponent $\alpha$.
In particular, $C^{0,1}$ denotes the class of functions
satisfying the usual Lipschitz condition.

\begin{e-theo}
\label{removal}
Let $H:\D^2 \to (M,J')$ be a
$C^\infty$ smooth map from the bidisc with coordinates $(z,w)$
to a smooth almost complex manifold $M$ of complex dimension 2.
Let $\Sigma$ be the set of all critical points of $H$.
Let $J= H^*J'$ be the pull-back of $J'$ on $\D^2 \backslash \Sigma$.
Suppose
\begin{itemize}
\item[(i)] for every $z \in \D$, the map
$\D\ni w \mapsto H(z,w)\in M$ is a $(J_{\st},J')$-holomorphic
immersion;
\item[(ii)] for every $z \in \D$, the set
$\{ z \} \times\D$ is not contained in $\Sigma$;
\item[(iii)] the map $H|_{\D^2\setminus\Sigma}$ preserves
the canonical orientations defined by $J_{st}$ and $J'$ on
$\D^2$ and $M$ respectively.
\end{itemize}
Then for every $z \in \D$, the set
$\Sigma\cap(\{z\}\times\D)$ is discrete, and
the almost complex structure $J$
defined on $\D^2 \backslash \Sigma$ extends to be continuous
on the whole bidisc $\D^2$.
Moreover, on every compact $K\subset\D^2$,
for some $0 < \alpha <1$,
the extension of $J$ is $C^\alpha$ in $z$ uniformply in $w$
and $C^{0,1}$ (Lipschitz) in $w$ uniformly in $z$.
If the structure $J'$ is integrable, then the extension of
$J$ is $C^\infty$ smooth on $\D^2$.
\end{e-theo}

{\bf Remarks.}

1. For simplicity we assume that all objects in the hypotheses
of Theorem \ref{removal} are smooth of class $C^\infty$,
however, the proof goes through for finite smoothness.
We leave the details to the reader.
The theorem is used in applications for constructing
pseudo-holomorphic discs in convenient coordinates.
After returning to the original manifold, the resulting
discs will be automatically smooth of class $C^\infty$
due to ellipticity.

2. In some applications we use a version of Theorem \ref{removal}
in which the map $H$ is smooth on $\bar\D\times\D$ and the
conclusion is that $J$ extends to all of $\bar\D\times\D$
with the stated regularity. That version formally does not follow
from Theorem \ref{removal}, but the proof goes through.

3. The condition (iii) can be replaced by (iii'):
{\it the set $\D^2\setminus\Sigma$ is connected.}
Indeed, if (iii') holds but (iii) does not, then it means
that $H$ changes the orientation to the opposite.
Let $\sigma(z,w)=(\bar z, w)$.
Then $H\circ\sigma$ satisfies (i--iii), and the conclusion of
Theorem \ref{removal} holds for $H\circ\sigma$,
whence for $H$.
The conditions (ii) and (iii) can be replaced
without much loss by a single condition
(ii'): {\it for every $z \in \D$, the set
$\Sigma\cap(\{z\}\times\D)$ is discrete.}
On the other hand, without (iii), Theorem \ref{removal}
fails. We give an example to that effect in
Section 3.2.

4. The resulting structure $J$ in Theorem \ref{removal}
does not have to be smooth in $w$.
We include an example to this effect in Section 3.2.
However, we do not know whether the smoothness
in $z$ can drop below $C^{0,1}$.
We admit that our proof of regularity
of $J$ in $z$ might not fully use all the hypotheses
of the theorem. Fortunately, the H\"older continuity
of $J$ in $z$ suffices for our applications.

5. Finally, if $J'$ is integrable, then a version of
Theorem \ref{removal} holds in higher dimension.
In that version, if $\dim_{\cc}M=n\ge2$, then
$z\in\D$ and $w\in\D^{n-1}$.
We leave the details to the reader.

\subsection{Reduction to PDE}

The condition (ii) in Theorem \ref{removal} has global
character. Nevertheless, we observe that it suffices to prove
Theorem \ref{removal} locally. More precisely, for $z\in\D$ put
$$
\Sigma_z=\{w\in\D:(z,w)\in\Sigma\}\qquad
K_z=\bar{\D\setminus\Sigma_z}\cap\D.
$$
According to this definition, if $w\in K_z$,
then (ii) holds in every neighborhood of $(z,w)$.
\begin{e-lemme}
\label{localization}
Suppose the conclusions of Theorem \ref{removal}
hold in a neighborhood of every point $(z,w)\in\D^2$
such that $w\in K_z$. Then they hold in all of $\D^2$.
\end{e-lemme}
\proof
Fix $z\in\D$.
By definition $K_z$ is closed in $\D$.
But it is also open because Theorem \ref{removal}
concludes that $\Sigma_z$ is discrete.
Hence $K_z=\D$. Since $z\in\D$ is arbitrary,
then Theorem \ref{removal} holds in all of $\D^2$.
\medskip

In the following proposition for simplicity
we add more assumptions to Theorem \ref{removal}.
In the proof of Theorem \ref{removal} we will use
this result locally.

\begin{e-pro}
\label{reductionPDE}
In addition to the hypotheses of Theorem \ref{removal}
we assume that $H$ is smooth on $\bar\D^2$ and
for every $z\in\bar\D$,
the map $w\mapsto H(z,w)$ is an embedding on $\bar\D$.
Then the complex matrix $A$ of $J$ on $\D^2\backslash\Sigma$
has the form
\begin{eqnarray}
\label{Aab}
 A = \left(
\begin{array}{cll}
a &  0\\
b &  0
\end{array}
\right),
\end{eqnarray}
where $a=g/f$, $b=ah_1+h_2$ for some
$f,g,h_1,h_2\in C^\infty(\D^2)$
satisfying the inequality $|f|\ge|g|$.
The singular set $\Sigma$ has the form
$\Sigma=\{|f|=|g|\}$, and for some
$\mu\in C^\infty(\D^2)$, the following system\ holds:
\begin{eqnarray}
\label{PDEfg}
f_{\bar w}=\mu\bar g, \qquad
g_{\bar w}=\mu\bar f.
\end{eqnarray}
\end{e-pro}
\proof
The statement involves $z$ as a parameter.
We first prove it for fixed $z$;
then it will be clear that the construction depends
smoothly on the parameter $z$ (see remark after the proof).

For simplicity put $z = 0$.
We introduce local coordinates $(z',w')$
in a neighborhood of the $J'$-complex curve $H(\{0\}\times\D)$
and use $(z'(z,w),w'(z,w))$
for the coordinate representation of $H$.
We choose the coordinates $(z',w')$ so that
\begin{eqnarray}
\label{PDEfg1}
z'(0,w) = 0,\qquad
w'(0,w) = w,
\end{eqnarray}
and for every $w'\in\D$,
the map $z'\mapsto(z',w')$ is $(J_\st,J')$-holomorphic.
Then the coordinate system $(z',w')$ preserves the
orientation of $M$ defined by $J'$.
Furthermore, $J'(0,w')=J_\st$, and the complex matrix $A'$ of
$J'$ satisfies
\begin{eqnarray}
\label{PDEfg2}  A'(0,w') = 0.
\end{eqnarray}
Using \cite{SuTu} (Lemma 2.4),
we modify the coordinates $(z',w')$ so that in addition to
(\ref{PDEfg1}) and (\ref{PDEfg2}) we have
\begin{eqnarray}
\label{PDEfg3}
 A'_{z'}(0,w') = 0.
 \end{eqnarray}
Put $Z=(z,w)$, $Z'=(z',w')$.
Then the complex matrix $A$ of the pull-back structure $J$
is obtained by the following transformation rule
(\cite{SuTu}, Lemma 2.3):
\begin{eqnarray}
\label{PDEfg4} A = (Z_Z' - A' \overline{Z}_Z')^{-1} (A'
\overline{Z}_{\overline Z}' - Z_{\overline Z}')
\end{eqnarray}
whenever this formula makes sense.
We want to describe $A(0,w)$.
By (\ref{PDEfg1}), (\ref{PDEfg2}) and (\ref{PDEfg4}),
$A(0,w)$ has the form
$A =-(Z_Z')^{-1}Z_{\overline Z}'$.
By (\ref{PDEfg1}) we have
\begin{eqnarray*}
Z'_Z =\left(
\begin{array}{cll}
z'_z &  0\\
w'_z &  1
\end{array}
\right) \qquad
Z'_{\bar Z} =\left(
\begin{array}{cll}
z'_{\bar z} &  0\\
w'_{\bar z} &  0
\end{array}
\right)
\end{eqnarray*}
We denote (for fixed $z=0$)
\begin{eqnarray}
\label{fgh}
f=z'_z,\;\;
g=-z'_{\bar z},\;\;
h_1=-w'_z,\;\;
h_2=-w'_{\bar z}.
\end{eqnarray}
The real Jacobian of the map $Z\mapsto Z'$
has the form $|f|^2-|g|^2$, hence by (iii)
we have $|f|\ge|g|$ and $\Sigma=\{|f|=|g|\}$.
Then $f\ne0$ on $\D^2\setminus\Sigma$,
and we immediately obtain the form (\ref{Aab})
of the matrix $A$ with expressions for $a$ and $b$.

We now derive the differential equations (\ref{PDEfg})
for $f$ and $g$.
The condition (i) of Theorem \ref{removal} in our coordinates
takes the form
$$
\left(
\begin{array}{cl}
z'\\
w'
\end{array}
\right)_{\overline w} = A'(z',w')
\left(
\begin{array}{cl}
\overline{z}'\\
\overline{w}'
\end{array}
\right)_{\overline w}.
$$
Differentiating this equation with respect to $z$,
since $A'(0,w') = 0$, we obtain for $z=0$
\begin{eqnarray*}
\left(
\begin{array}{cl}
z'\\
w'
\end{array}
\right)_{\overline{w}z} = A_z'\left(
\begin{array}{cl}
0\\
1
\end{array}
\right).
\end{eqnarray*}
We have
$A'_z =A'_{z'}z'_z + A'_{\overline z'} \overline{z}'_z
+ A'_{w'}w'_z + A'_{\overline w'} \overline{w}'_z
=A'_{\overline z'} \overline{z}'_z$ because
$A'_{z'}(0,w') = 0$ by (\ref{PDEfg3}), and
$A'_{w'}(0,w') =A'_{\overline w'}(0,w') = 0$ by
(\ref{PDEfg2}).
Hence
\begin{eqnarray*}
\left(
\begin{array}{cl}
z'\\
w'
\end{array}
\right)_{\overline w z} = A'_{\overline z'} \left(
\begin{array}{cl}
0\\
1
\end{array}
\right) \overline{z}'_z.
\end{eqnarray*}
Let $\mu$ denote the $(1,2)$ entry of the matrix
$-A'_{\overline{z}'}(0,w)$. Then
$z'_{z\overline w}(0,w) = -\mu \overline{z}'_z(0,w)$.
Using the notation (\ref{fgh}), we immediately obtain the
first equation in (\ref{PDEfg}). The second equation
in (\ref{PDEfg}) is derived similarly.
It remains to add that our construction including
the choice of the coordinates $Z'$ depends smoothly
on the parameter $z$.
Proposition \ref{reductionPDE} is proved.
\medskip

{\bf Remarks.}

1. In the above proof, we use a version of Lemma 2.4
from \cite{SuTu} with smooth dependence on parameters.
A careful examination of the proof in \cite{SuTu} shows that
the desired version holds. In particular, we recall that the only
analytic tool used in the proof is solving the equation
$u_{\bar w}=p(w)u+q(w)$.
This is similar to solving an ordinary differential equation
$dy/dx=p(x)y+q(x)$. In the procedure of solving this equation,
one replaces integration with
respect to $x$ by the Cauchy-Green integral (\ref{CauchyGreen}).
The latter is known to depend smoothly on parameters.
Hence if the coefficients $p$ and $q$ smoothly depend
on additional parameters, then there exists a solution $u$
that smoothly depends on the parameters.

2. We can now conclude the proof of Theorem \ref{removal}
in the important special case, in which the structure $J'$
is integrable.
By Lemma \ref{localization} it suffices to prove the result locally.
Then we can use Proposition \ref{reductionPDE}.
In its proof we have $A'=0$, hence $\mu=0$, and the functions
$f$, $g$ are holomorphic in $w$.
Then, by the maximum principle, one can see that
$\Sigma=\{f=0\}$, and
$\Sigma\cap(\{ z \} \times\D)$ is discrete.
By the removable singularity theorem,
the ratio $a=g/f$ is holomorphic in $w$
on the whole bidisc.
Then in fact $a$ is $C^\infty$ smooth in both $z$ and $w$
by the Cauchy integral formula in $w$.
By the maximum principle, $|a|<1$ holds for the extension.
By Proposition \ref{JtoA}, the matrix (\ref{Aab}) defines
an almost complex structure if and only if $|a|\ne1$.
Hence the extension of $J$ is well defined and $C^\infty$,
which concludes the proof. In Appendix we include a proof
for integrable structures independent of
Proposition \ref{reductionPDE}.

\subsection{Two examples}

The following simple example shows that the condition (iii)
in Theorem \ref{removal} cannot be omitted.
\medskip

{\bf Example.}
Let $M=\C^2$, $J'=J_{st}$.
Define $H:\D^2\to M$ by
$$
z'=z-2\bar zw, \qquad
w'=w.
$$
Then $f=z'_z=1$, $g=-z'_{\bar z}=2w$, $a=g/f=2w$.
The real Jacobian of $H$ has the form
$|f|^2-|g|^2=1-4|w|^2$. It vanishes on the
real hypersurface $\Sigma=\{|w|=1/2\}$.
Then $H^*J'$ can not be extended to $\Sigma$ because
$|a|=1$ on $\Sigma$. The conditions (iii) and (iii')
are not fulfilled.
\medskip

The following example shows that the drop of
smoothness with respect to $w$ in Theorem \ref{removal}
can occur.
\medskip

{\bf Example.}
Let $M=\C^2$ with coordinates $(z',w')$. Let the
almost complex structure $J'$ have the complex matrix
$$
A' = \left(
\begin{array}{cll}
\bar w' &  -\bar z'\\
0 &  0
\end{array}
\right).
$$
Consider a blow-up map
$Z=(z,w)\mapsto Z'=H(z,w)=(zw,w)$.
We find the complex matrix $A$ of the pull-back
$J=H^*J'$ by (\ref{PDEfg4}).
Since $H$ is holomorphic in the usual sense,
we have $A=(Z'_Z)^{-1}A'\bar{Z'_Z}$, which yields
$$
A = \left(
\begin{array}{cll}
w^{-1}\bar w^2 &  0\\
0 &  0
\end{array}
\right).
$$
The map $H$ satisfies the hypotheses of Theorem \ref{removal}.
In particular, for fixed $z$ the map $w\mapsto H(z,w)$
is $J'$-holomorphic because the matrix $A$ has zeros
in the second column. The singular set
$\Sigma$ is the line $w=0$. We realize that $A$, whence $J$
is not smooth but merely Lipschitz in $w$ in accordance
with Theorem \ref{removal}.

\subsection{H\"older continuity of the logarithmic difference}

We consider the equation
\begin{eqnarray}
\label{Vekua}
h_{\bar w}=\mu\bar h
\end{eqnarray}
in a bounded domain $G\subset \cc$.
Although in our applications $\mu$ will be smooth,
one can assume that $\mu$ is merely bounded and
$h_{\bar w}$ in (\ref{Vekua}) is a Sobolev derivative.
The equation is relevant because both $f+g$ and $f-g$
for $f$ and $g$ in (\ref{PDEfg}) satisfy an equation
of the form (\ref{Vekua}), which we will use later.
Solutions of (\ref{Vekua}) are called
{\it generalized analytic functions} in \cite{Ve}.
They have the following representation
\begin{eqnarray}
\label{VekuaRep}
h=\phi e^{Tu}, \qquad
u=\mu\bar h/h.
\end{eqnarray}
Here $T = T_G$ denotes the Cauchy--Green integral
\begin{eqnarray}
\label{CauchyGreen}
T u(w) = \frac{1}{2\pi i} \int\int_\D
\frac{u(\tau)\,d\tau \wedge d\overline\tau}{\tau -w}.
\end{eqnarray}.

The function $\phi$ is holomorphic in $G$.
Indeed, since $\partial_{\bar w}Tu=u$, then
$$
\partial_{\bar w}\phi=\partial_{\bar w}(he^{-Tu})
=\mu \bar he^{-Tu}+he^{-Tu}(-u)=0.
$$
In particular, the zero set of $h$ is discrete
unless $h\equiv0$.
The function $Tu$ is called the
{\it logarithmic difference} of $h$ because it measures
the distance from $h$ to a holomorphic function $\phi$
in the logarithmic scale.

Since $\mu$ is bounded, then the logarithmic difference
of $h$ and $h$ itself are bounded in the H\"older norm in $w$.
We now obtain the following result about H\"older
continuity of the logarithmic difference on a parameter.

\begin{e-pro}
\label{Holder}
In the closed bidisc $\overline\D^2$ with coordinates
$(z,w)$, let $h,\mu \in C^\infty(\overline\D^2)$ satisfy
(\ref{Vekua}).
Suppose $h \neq 0$ on $\{(z,w): |w|=1\}$.
Then $h$ has the representation (\ref{VekuaRep}),
in which $\phi \in C^\infty(\overline\D^2)$
and holomorphic in $w$.
Furthermore,
$T u \in C^\alpha(\overline\D^2)$
for some $0<\alpha<1$.
(The operator $T=T_{\D}$ is applied
with respect to $w$.)
\end{e-pro}

{\bf Remark.}
In the proof we will obtain an estimate
$\alpha= 1/(n+1)$, where $n$ locally is the maximum number
of zeros of $h$ in $w$. We do not know whether this
estimate is sharp.
\medskip

We need two lemmas in the proof.
We use the notation $d^2w=\frac{i}{2}dw\wedge d\bar w$
for the area element.
We denote by $\m(E)$ the area of $E\subset\cc$.

\begin{e-lemme}
\label{rotation} For every measurable set $E\subset\cc$,
we have
$\int\int_E |w|^{-1}d^2w \leq 2(\pi\,\m(E))^{1/2}$.
\end{e-lemme}

\proof
We have
$I =\int\int_E |w|^{-1}d^2w \leq
\int\int_{|w|< r}|w|^{-1}d^2w$,
where $\m(E) = \pi r^2$.
Then by evaluating the last integral explicitly
and expressing $r$ in terms of $\m(E)$, we get
$I \leq 2\pi r = 2(\pi\,\m(E))^{1/2}$
as desired.
\bigskip

\begin{e-lemme}
\label{area}
Let $p(w) = (w - w_1)\ldots(w-w_n)$, and let
$A(\delta)  =\m\{ w: \vert p(w) \vert < \delta \}$.
Then $A(\delta) \leq \pi n\delta^{2/n}$.
\end{e-lemme}

\proof Let $\vert p(w) \vert < \delta$. Then $\vert w - w_j \vert
< \delta^{1/n}$ for some $j$.
Then $w \in \cup_k \{w: |w - w_k|< \delta^{1/n} \}$,
and the lemma follows.
\bigskip

{\bf Proof of Proposition \ref{Holder}:}
We use the notation $C_1, C_2, \dots$ for constants.
We have
$u = \mu\overline h /h = v/h$, where $v = \mu \overline h$.
Then $|u| \leq C_1 = ||\mu||_\infty$.
Since $u$ is bounded, then obviously
$Tu(z,\bullet)\in C^\alpha$
for every $0<\alpha<1$ uniformly in $z$.

We need to prove that $Tu$ is $C^\alpha$ in $z$ for
some $0<\alpha<1$ uniformly in $w$.
Set $\Delta z = z' - z''$.
Omitting $w$ for simplicity, we have
\begin{eqnarray*}
& & |u(z') - u(z'')| = |(v/h)(z') - (v/h)(z'')|
= \left|\frac{v(z') - v(z'')}{h(z')} -
\frac{v(z'')}{h(z'')}\cdot\frac{h(z') - h(z'')}{h(z')} \right|\\
& &\leq \frac{C_2|\Delta z|}{|h(z')|} + C_1
\frac{C_2|\Delta z|}{|h(z')|}
= C_3\frac{|\Delta z|}{|h(z')|}.
\end{eqnarray*}

Set $\Delta T = \vert Tu(z',w_0) - Tu(z'',w_0) \vert$.
Using Lemma \ref{rotation} for the second integral below,
\begin{eqnarray*}
& &\Delta T \leq C_3 |\Delta z|\, \delta^{-1}
\int\int_{|h(z',w)| > \delta, |w| < 1}
\frac{d^2w}{|w - w_0|}
+ 2C_1 \int\int_{|h(z',w)| < \delta, |w|<1}
\frac{d^2w}{|w - w_0|}\\
& &\leq C_4\left(|\Delta z|\, \delta^{-1}
+ \m\{ w: |h(z',w)| < \delta, |w|<1 \}^{1/2}
\right ).
\end{eqnarray*}

Define $\phi=h e^{-Tu}$.
Then $\phi$ is holomorphic in $w$.
Note that $Tu$ is $C^\infty$ in $(z,w)$
outside the zero set of $h$, in particular,
for $|w|=1$. By the Cauchy
integral formula, $\phi\in C^\infty(\overline\D^2)$.

Since $u$ is bounded, we put $|Tu| \leq C_5$.
Then
$$
|h| e^{-C_5} \le |\phi| \le |h|e^{C_5}.
$$
Let $w_1,\ldots,w_n$ be the zeros of
$h(z,w)$ in $\D$ for fixed $z$, and let
$p(w) = (w -w_1)\ldots(w-w_n)$.
By the argument principle for $\phi$, the number
$n$ does not depend on $z$.
For $|w|=1$ we have $|h|\ge C_6 > 0$,
and $|p|\le 2^n$.
Then by the minimum principle
\begin{eqnarray*}
|\phi p^{-1}| \ge
C_6e^{-C_5}2^{-n} = C_7>0.
\end{eqnarray*}

The condition $|h|\le \delta$ implies
$|p|\le|\phi|/C_7 \le C_8 |h| \le C_8 \delta$.
By Lemma \ref{area} we have the estimate
$\m\{|h| \le \delta \} \le C_9 n \delta^{2/n}$ and
$\Delta T\le C_{10} \left(|\Delta z|\,\delta^{-1}
+\delta^{1/n} \right )$.
Put
$\delta = |\Delta z|^{1-\alpha}$. Then
$\Delta T \leq C_{11} \left(|\Delta z|^{\alpha}
+|\Delta z|^{(1-\alpha)/n} \right )$.
Take now $\alpha = 1/(n+1)$. Then
$\Delta T \leq C_{12}|\Delta z|^{\alpha}$.
Thus $T u \in C^\alpha(\overline\D^2)$,
which concludes the proof of Proposition \ref{Holder}.

\subsection{Decomposition of the phase of a complex polynomial}

The results of Sections 3.4--3.5 are needed only for the proof
that $A$, hence $J$, is Lipschitz in $w$.
This is used in the proofs of the results of Section 5
about gluing discs to real tori, which are not immersed in general.

We call $\<w \> := {\overline w}/w$ the {\it phase function}
of $w\in\cc$. Let
$$\Delta_n = \{ t = (t_1,...,t_n): t_j \geq 0, \sum t_j = 1 \}$$
be the standard $(n-1)$ simplex.

\begin{e-theo}
\label{interpolation1}
For every integer $n \ge 1$ there exists a
constant $C_n>0$ and measures $\mu_{nk}$, $1\le k\le n$, on
$\Delta_n$ depending on parameters $w_1,...,w_n \in \cc$ such that
$\int_{t \in \Delta_n} \vert d\mu_{nk}\vert \leq C_n$ and the
following identity holds:
\begin{eqnarray*}
\<(w - w_1)\dots(w - w_n)\> = \sum_{k=1}^n \int_{t \in\Delta_n}
\< w - t_1 w_1 -\dots- t_n w_n \>^k
d\mu_{nk}(w_1,\dots,w_n,t)
\end{eqnarray*}
\end{e-theo}

The above formula can be made much more precise. The
singular measures $\mu_{nk}$ reduce to integration over some
subsimplexes of $\Delta_n$ with bounded densities.
Theorem \ref{interpolation1} means that the phase function
of a polynomial can be decomposed into a ``sum'' of
the phase functions of binomials.
It is somewhat similar to decomposition of rational
functions into partial fractions, but the sum
in fact turns into an integral.
We first prove the result in a special case, in which
the polynomial is a product of just two binomials.
Then the general case will follow by induction.

\begin{e-lemme}
\label{interpolation2}
For every integer $n\ge1$, $1\le k\le n$, $1\le j\le n+1$,
there are constants $c_{kj}\in\R$, such that for every
$w,w_0 \in \cc$, the following identity holds
\begin{eqnarray*}
\<w^n(w - w_0)\>
= \<w_0 w^n\> +\<w_0^n (w -w_0)\>
+ \sum_{k=1}^n\sum_{j=1}^{n+1} c_{kj}
\<w_0\>^{n+1-j}\int_0^1 \<w - w_0t\>^j
(1-t)^{k-1}dt.
\end{eqnarray*}
\end{e-lemme}

\proof
Put
$\Phi = \<w^n(w-w_0)\>$.
We will use the partial fraction decomposition
\begin{eqnarray*}
\frac{1}{w^n(w-w_0)} = \frac{1}{w_0^n(w-w_0)} - \frac{1}{w_0^n w}
-\dots- \frac{1}{w_0 w^n}.
\end{eqnarray*}
Then
\begin{eqnarray*}
&&\Phi= \frac{\overline{w}^n(\overline{w} -
\overline{w}_0)}{w_0^n(w-w_0)} - \sum_{k=1}^n
\frac{\overline{w}^n(\overline{w} -
\overline{w}_0)}{w_0^{n-k+1}w^k} = \frac{(\overline{w}_0 +
(\overline{w} - \overline{w}_0))^n(\overline{w} -
\overline{w}_0)}{w_0^n(w-w_0)}- \sum_{k=1}^n
\frac{\overline{w}^n(\overline{w} -
\overline{w}_0)}{w_0^{n-k+1}w^k}\\
&&= \sum_{k=0}^n
\left(
\begin{array}{cl}
n\\
k
\end{array}
\right)
\frac{\overline{w}_0^{n-k}(\overline{w} -
\overline{w}_0)^{k+1}}{w_0^n(w-w_0)} -\sum_{k=1}^n
\frac{\overline{w}^n(\overline{w} -
\overline{w}_0)}{w_0^{n-k+1}w^k}.
\end{eqnarray*}
Put
$$
A^k_l(w)= \frac{\overline{w}_0^{n-k+1}
\overline{w}^k}{w_0^{n-l+1}w^l}.
$$
Then
\begin{eqnarray*}
\Phi = A_1^1(w - w_0) + A_n^n(w) + \sum_{k=1}^n \left(
\begin{array}{cl}
n\\
k
\end{array}
\right) A_1^{k+1}(w-w_0)
-\sum_{k=1}^n A_k^{n+1}(w)
+\sum_{k=1}^{n-1}A_k^n(w).
\end{eqnarray*}
The terms $A_1^1(w-w_0)$ and $A_n^n(w)$ are listed separately
because they are the only bounded terms in the above formula. All
other terms have the form $A^p_q$ with $p > q$ so they are
unbounded as $w_0\to0$.
Put
$$
f(t) = A_1^{k+1}(w - t w_0).
$$
We use Taylor's expansion
\begin{eqnarray*}
f(1) = \sum_{n=0}^{k-1} \frac{1}{n!}f^{(n)}(0) + \frac{1}{(k-1)!}
\int_0^1 f^{(k)}(t)(1-t)^{k-1}dt.
\end{eqnarray*}
to express $A_1^{k+1}(w-w_0)$ in terms of $A^p_q(w)$.
Note that
\begin{eqnarray*}
& &\frac{d}{dt} A^p_q(w - tw_0) = q A^p_{q+1}(w - t w_0) - p
A^{p-1}_q(w - t w_0),\\
& &\frac{d}{dt}\bigg\vert_{t=0} A^p_q(w - t w_0)
= q A^p_{q+1}(w) - pA^{p-1}_q(w).
\end{eqnarray*}
By Taylor's formula, not keeping track of the exact values of the
coefficients, we have
\begin{eqnarray*}
A_1^{k+1}(w - w_0) = \sum_{1 \leq q < p \leq
k+1}\alpha_{kq}^pA^p_q(w) + \sum_{p=1}^{k+1}\beta_{kp} \int_0^1
A^p_p (w - t w_0) (1 - t)^{k-1}dt
\end{eqnarray*}
where $\alpha^p_{kq}$ and $\beta_{kp}$ are universal constants. Then

\begin{eqnarray*}
\Phi = A_1^1(w - w_0) + A^n_n(w) + \sum_{1 \leq q < p \leq n+1}
a^p_qA^p_q(w) + \sum_{k=1}^n\sum_{p=1}^{n+1}
c_{kp} \int_0^1 A^p_p(w - t w_0)(1 - t)^{k-1}dt
\end{eqnarray*}
where $a^p_q$ and $c_{kp}$ are universal constants.

Since $|A^p_p|= 1$ and $|\Phi| = 1$, the
sum $\sum_{1 \leq q < p \leq n+1} a^p_q A^p_q(w)$ is
bounded.  But the terms $A^p_q(w)$ with $q<p$
are all unbounded and have different asymptotics as $w_0\to0$.
Hence $a^p_q = 0$, and Lemma \ref{interpolation2}
follows.

\subsection{Lipschitz continuity of
the logarithmic difference}

\begin{e-pro}
\label{final} For every $\mu \in C^1(\overline\D)$,
$\varepsilon >0$, $M > 0$ and integer $n\ge0$
there exists a constant $C >0$ such that every function
$h \in C^1(\overline\D)$ satisfying the conditions
\begin{itemize}
\item[(i)] $h_{\overline w} = \mu \overline h$
\item[(ii)] $|h(w)|> \varepsilon$ for $|w| > 1/2$
\item[(iii)] $h$ has $n$ zeros in $\D$
\item[(iv)] $||h||_{C^1(\overline\D)} \leq M$
\end{itemize}
admits the representation
$h = \phi_0 p e^{T u}$, where
$u =\mu\overline h/h$,
$p$ is a monic polynomial of degree $n$,
and we have the estimates
$||\phi_0||_{C^1(\overline\D)} \leq C$,
$|\phi_0|\geq 1/C$, $||T u||_{C^{0,1}(\D)} \leq C$.
\end{e-pro}

The goal of Proposition \ref{final} is that
the estimates on $\phi_0$ and $Tu$ depend only on
the number of zeros of $h$, not their location.

\begin{e-lemme}
\label{Holder1}
Let $\lambda \in C^{\alpha}(\cc)$ for some
$0 <\alpha < 1$ and $\lambda(0) = 0$.
Then for every positive integer $n$ we have
$$||\lambda \langle w \rangle^n||_{C^\alpha(\cc)}
\leq n C||\lambda||_{C^\alpha(\cc)}$$
where $C> 0$ is an absolute constant.
\end{e-lemme}

Lemma \ref{Holder1} follows from a more general result
\cite{SuTu} (Lemma 5.4), which in place of $\<w\>^n$
has a function whose derivatives have the estimate
$O(|w|^{-1})$.

\begin{e-lemme}
\label{Cauchy}
For every positive integer $n$ and $w,w_0 \in \D$ we have
\begin{eqnarray}
\label{Cauchy1}
T \< w - w_0 \>^n = \frac{1}{n+1}
\frac{(\overline{w}-\overline{w}_0)^{n+1}}{(w - w_0)^n}
\end{eqnarray}
\end{e-lemme}
\proof
Let $f(w) = \langle w - w_0 \rangle^n$, and let
$g(w)$ be the right-hand part of (\ref{Cauchy1}).
Then $g_{\overline w}=f$.
By the Cauchy-Green formula,
$g(w) = K g (w) + T f(w)$, $w \in \D$.
Here
$$
K g (w) =\frac{1}{2 \pi i} \int_{|\zeta|=1}
\frac{g(\zeta)d\zeta}{\zeta - w}
$$
is the Cauchy type integral over the
unit circle.
But for $|w|=1$ we have
$$
g(w) = \frac{1}{n+1}(w^{-1} - \overline{w}_0)^{n+1} \left[
w^{-1} \sum_{k=0}^\infty \left ( \frac{w_0}{w} \right )^k
\right]^n.
$$
Thus the Laurent series of $g$ on the unit circle contains
only negative powers of $w$. Hence $Kg \equiv 0$,
and the lemma follows.
\bigskip

\begin{e-lemme}
\label{Holder2}
Let $u =\lambda \<p \>$, where $p(w) = (w-w_1)...(w-w_n)$,
$\lambda \in C^\alpha(\D)$. Then
$$
||T u||_{C^{0,1}(\D)} \le C||\lambda||_{C^\alpha(\D)},
$$
where $C$ depends on $n$ and $\alpha$ only.
\end{e-lemme}
\proof
We represent $\<p \>$ by
Theorem \ref{interpolation1}.
(This is the only instance when Theorem \ref{interpolation1}
is used.)
Since the integrals
$\int |d\mu_{nk}|$ are uniformly bounded, it suffices to prove the
result separately for each term in Theorem \ref{interpolation1}.
Hence it suffices to consider the case where
$$
p(w) = (w - w_0)^k,\quad
|w_0|< 1, \quad
1 \leq k \leq n.
$$
Then
$u(w)=(\lambda(w)-\lambda(w_0))\<w-w_0\>^k
+\lambda(w_0)\<w-w_0\>^k$, hence
$$
Tu(w)=T[(\lambda(w)-\lambda(w_0))\<w-w_0\>^k]
+\frac{\lambda(w_0)}{k+1}
\frac{(\overline{w}-\overline{w}_0)^{k+1}}{(w-w_0)^k}.
$$
The first term is uniformly bounded in $C^{1,\alpha}(\D)$,
because $(\lambda(w)-\lambda(w_0))\<w-w_0\>^k$
is uniformly bounded in $C^\alpha(\D)$ by Lemma \ref{Holder1}.
The second term is obtained by Lemma \ref{Cauchy};
clearly, it is in $C^{0,1}(\D)$. This proves the lemma.
\medskip

{\bf Proof of Proposition \ref{final} :}
Without loss of generality
$||\mu||_{C^1(\overline\D)} \leq M$.
We use $C_1, C_2,\dots$ for constants
depending on $\varepsilon$, $M$ and $n$ only.
In this proof, the term ``uniformly bounded'' means
bounded by a constant depending on
$\varepsilon$, $M$ and $n$ only.

Since $h$ satisfies the equation (i),
then it admits the representation
$h=\phi e^{Tu}$, $u=\mu\overline h/h$ with holomorphic $\phi$.
Set $\phi_0=\phi/p=he^{-T u}/p$, where
$p(w) = (w - w_1)\dots(w- w_n)$, $w_j$
are the zeros of $h$, $|w_j|< 1/2$.
Then $u = \lambda \<p\>$, where
$\lambda=(\mu\overline{\phi}_0/\phi_0) e^{\overline{T u} - T u}$.
We will see that
$\lambda$ is uniformly bounded in $C^\alpha(\overline\D)$,
and Proposition \ref{final} will follow by Lemma \ref{Holder2}.

Fix any $0 < \alpha < 1$, say, $\alpha = 1/2$.
Since $|\mu|\le M$, then $|u|\leq M$.
Since the operator $T: L^\infty(\D) \to C^\alpha(\D)$
is bounded, we have
$||T u||_{C^\alpha(\D)} \le C_1$ and
$||e^{-T u}||_{C^\alpha(\D)} \le C_2$.
The condition $\varepsilon \le|h(w)|\le M$
for $w\in b\D$ implies the inequality
$1/C_3 \le|\phi(w)|\le C_3$ for $w \in b\D$.
Since $|w_j|< 1/2$, then
$2^{-n} \le|p(w)|\le 2^n$ for $w \in b\D$.
Therefore for
$\phi_0 = \phi/p$ we obtain $1/C_4 \le|\phi_0(w)|\le C_4$
for $w \in b\D$, hence for all $w \in \overline\D$ because
$\phi_0$ is holomorphic and has no zeros.
We now show that $\parallel \phi_0 \parallel_{C^1(\overline\D)}$
is uniformly bounded. Since
$$
(p^{-1})' = -p^{-1}\sum_{j=1}^n
(w-w_j)^{-1},
$$ then $||p^{-1}||_{C^1(b\D)} \leq C_5$.
By splitting
$\int\int_\D =\int\int_{(1/2)\D}
+\int\int_{\D \backslash(1/2)\D}$ we obtain
$$
||T u||_{C^1(b\D)} \leq C_6 (||u||_{L^\infty((1/2)\D)}
+||u||_{C^1(\overline{\D} \backslash (1/2)\D}).
$$
The last term has the estimate
$$
||u||_{C^1(\bar\D \backslash (1/2)\D} \leq M
||\overline{h}/h||_{C^1(\bar\D \backslash (1/2)\D)}
\le C_7 \varepsilon^{-1} ||h||_{C^1(\bar\D)}
\le C_8.
$$
Therefore
$||T u||_{C^1(b\D)} \le C_9$ and
$||e^{-T u}||_{C^1(b\D)} \leq C_{10}$.
Then for
$\phi_0 =e^{-Tu}h/p$ we obtain
$||\phi_0||_{C^1(\overline\D)} \le C_{11}$.
Since $|\phi_0|> 1/C_4$ we have
$||\overline{\phi}_0/\phi_0||_{C^1(\overline\D)} \le C_{12}$.
Now for
$\lambda=(\mu\overline{\phi}_0/\phi_0) e^{\overline{T u}-T u}$,
we have
$||\lambda||_{C^\alpha(\overline\D)} \le C_{13}$,
as desired. Proposition \ref{final} is proved.

\subsection{Proof of Theorem \ref{removal}}

We resume the proof of Theorem \ref{removal} and return
to the notation of Section 3.1.
by Lemma \ref{localization} it suffices to prove
Theorem \ref{removal} locally.
Hence it suffices to prove it in the settings of
Proposition \ref{reductionPDE}.
Recall that $\Sigma=\{|f|=|g|\}$, where $f$ and $g$ satisfy
the equations (\ref{PDEfg}) and the inequality
$|f|\ge|g|$. Put $\Sigma'=\{f=0\}$. Then $\Sigma'\subset\Sigma$.

\begin{e-lemme}
\label{SigmaIsDiscrete}
For every $z \in \D$, the set
$\Sigma'\cap(\{z\}\times\D)$ is discrete.
\end{e-lemme}
\proof
We prove the lemma for fixed $z$ and treat
$\Sigma'$ as a subset of $\D$.
Let $G:=\D\setminus\Sigma'$.
Put $u=\mu\bar g/f$ in $G$.
Then $f_{\bar w}=uf$ in $G$.
Without loss of generality
(by the hypotheses of Proposition \ref{reductionPDE})
we assume that $\mu$ is bounded.
Since $|f|\ge|g|$, then $u$ is bounded in $G$.
Put $\phi=fe^{-Tu}$, where $T=T_G$ is the Cauchy-Green
integral (\ref{CauchyGreen}).
Then $\phi$ is continuous in $\D$ and holomorphic in $G$.
By the definition of $\phi$, we have
$\Sigma'=\{\phi=0\}$.
By Rado's theorem, $\phi$ is holomorphic on all of $\D$.
By (ii) of Theorem \ref{removal}, $\Sigma'\ne\D$.
Hence $\Sigma'$ is discrete.
The lemma is proved.

\begin{e-lemme}
\label{Sigma}
$\Sigma'=\Sigma$.
\end{e-lemme}
\proof
We again treat $\Sigma$ as a subset of $\D$.
Arguing by contradiction,
let $w_0\in\Sigma\setminus\Sigma'$.
Then $|f(w_0)|=|g(w_0)|\ne0$.
By multiplying $f$ by an appropriate constant of modulus $1$,
we can assume $f(w_0)=g(w_0)$.
Put $h = 1 - \frac{g}{f}$.
Then $\Re h\ge 0$.
Using (\ref{PDEfg}) we obtain that $h$ in $\D\setminus\Sigma'$
satisfies the following equation:
$$
h_{\overline w}=\frac{g f_{\bar w}-
f g_{\bar w}}{f^2} = \lambda h,\qquad
\lambda= -\mu{|f|\over f}
\left(1+{|g|\over|f|}\right)
{|f|-|g|\over f-g}.
$$
Note that $\lambda$ is bounded.
Then $h=\phi e^{T\lambda}$, where $\phi$ is holomorphic
in $\D\setminus\Sigma'$ and $\phi(0) = 0$.
In fact $\phi$ is holomorphic in all of $\D$
because it is bounded and has isolated singularities.
We claim that $h\equiv0$. Otherwise, since $T \lambda$
is continuous, then $h=\phi e^{T\lambda}$ maps
a neighborhood of $w_0$ onto a neighborhood of 0,
which is not possible because $\Re h\ge0$.
Now $h\equiv0$, that is, $f\equiv g$ contradicts
$\Sigma\ne\D$.
The lemma is proved.
\medskip

To complete the proof of Theorem \ref{removal}
we need to show that the function $a=g/f$
extends to all of $\D^2$ with the stated
regularity properties, and that the extension
satisfies $|a|<1$. Then by Proposition \ref{JtoA},
the complex matrix $A$ of the form (\ref{Aab}) will
define the desired extension of
$J=H^*J'$.

We first note that if $a$ extends continuously
to $\Sigma$, then $|a|<1$ follows immediately
by applying to $h=1-a$ the argument from
the proof of Lemma \ref{Sigma}.

We make the substitution
\begin{eqnarray}
\label{fgtilde}
\tilde f=f+g,\qquad
\tilde g=f-g,\qquad
\tilde a = \tilde g/\tilde f=\frac{1-a}{1+a}
\end{eqnarray}
and drop the tildes.
The new $f$, $g$ and $a$ are defined on $\D^2\setminus\Sigma$
and satisfy instead of (\ref{PDEfg}) the following equations
\begin{eqnarray}
\label{PDEfgtilde}
f_{\bar w}=\mu\bar f, \qquad
g_{\bar w}=-\mu\bar g. \qquad
\end{eqnarray}
We need to show that $a$ extends to all of $\D^2$
with the stated regularity.
It suffices to prove the theorem locally.
Thus without loss of generality we assume that
there exist $n$, $\varepsilon$, and $M$ such that
both $f$ and $g$ satisfy the hypotheses of Proposition \ref{final}.
Since $f$ and $g$ have the same zero set,
then they have the representations
$$
f= \phi_0 p e^{T u},\quad u = \mu \overline{ f}/f,
$$
$$
g = \psi_0 p e^{ T v},\quad v = - \mu \overline{ g}/g.
$$
Here $p(z,w) = (w - w_1(z))\dots(w - w_n(z))$, where
$w_1(z),\dots, w_n(z)$ are the zeros of $f$ and $g$
for fixed $z\in\D$.
Then by Propositions \ref{Holder} and \ref{final}, the formula
$a =g/f=(\phi_0/\psi_0) e^{ T u - T v}$
defines the extension of $a$ with the needed regularity.
The proof is complete.

\section{Gluing $J$-holomorphic discs to the standard torus}

Let $\D^2$ denote the standard bidisc in $\cc^2$
with coordinates $(z,w)$.
Let $J$ be an almost complex structure in $\D^2$
with complex matrix $A$ of the form (\ref{Aab}).
A map $\D\ni\zeta \mapsto (z(\zeta),w(\zeta))\in\D^2$
is $J$-holomorphic if and only if it satisfies
the following quasi-linear system:
\begin{eqnarray}
\label{mainsystem}
\left\{
\begin{array}{cccc}
& &z_{\overline\zeta}=a(z,w)\overline z_{\overline\zeta}\\
& &w_{\overline\zeta}=b(z,w)\overline z_{\overline\zeta}
\end{array}
\right.
\end{eqnarray}
We assume that $|a(z,w)|\le a_0<1$,
which implies the ellipticity of the system.
The following theorem strengthens one of the main
results of \cite{CoSuTu}.
For $r>0$ we put $\D_r:= r\D$.

\begin{e-theo}
\label{Riemann}
Let
$a,b:\overline\D\times\overline\D_{1+\gamma}\to\cc$,
$\gamma>0$. Let $0<\alpha<1$.
Suppose $a(z,w)$ and $b(z,w)$ are $C^\alpha$
in $z$ uniformly in $w$ and $C^{0,1}$ (Lipschitz)
in $w$ uniformly in $z$.
Suppose
\begin{eqnarray*}
|a(z,w)|\le a_0 < 1, \qquad
a(z,0)=0, \qquad
b(z,0)=0.
\end{eqnarray*}
Then there exist $C>0$ and integer $N\ge1$
such that for every integer $n \ge N$, real $0 < r \leq 1$ and
$0 \leq t < 2\pi$
(alternatively, there exist $C>0$ and $0<r_0\le1$
such that for every $n\ge0$, $0<r<r_0$ and $0 \leq t < 2\pi$),
the system (\ref{mainsystem}) has a unique solution
$(z_r,w_r):\overline\D\to\overline\D\times\overline\D_{1+\gamma}$
of class $C^{1,\alpha}$
with the properties:
\begin{itemize}
\item[(i)] $|z_r(\zeta)|=1$, $|w_r(\zeta)|=r$ for $|\zeta|= 1$;
$z_r(0) = 0$, $z_r(1)=1$ and $w_r(1)=re^{it}$;
\item[(ii)] $z_r: \overline \D \to \overline \D$
is a diffeomorphism;
\item[(iii)] $|w_r(\zeta)| \le Cr|\zeta|^n$,
and the winding number of $w_r|_{b\D}$ is equal to $n$;
\item[(iv)] for fixed $r$ we have
$\{(z_r(\zeta),w_r(\zeta)): |\zeta|=1,\, 0\le t<2\pi\}=b\D\times b\D_r$.
\end{itemize}
The solution continuously depends on the parameters $r$, $t$
and the coefficients $a$ and $b$.
In particular the map $ r  \mapsto (z_r,w_r)$ is a homotopy between
$(\zeta,0)$ and  $(z_1,w_1)$.
\end{e-theo}

We note that the conclusions (i--iii) remain true
if the coefficients $a$ and $b$ are in $C^\alpha$
in both $z$ and $w$, that is, without assuming
that they are $C^{0,1}$ in $w$.

The proof for smooth $a$ and $b$, $r = 1$, $n$ big enough
and no $t$ is given in \cite{CoSuTu}.
The proof of the present statement is similar;
we briefly describe it below.

We look for a solution of (\ref{mainsystem}) in the form
$z = \zeta e^u$, $w = re^{it}\zeta^n e^v$.
Then the new unknowns $u$ and $v$ satisfy
a similar system but with linear boundary conditions.

We reduce the system of PDE for $u$ and $v$ to a system of
singular integral equations using suitable modifications of
the Cauchy--Green operator and the Ahlfors--Beurling
transform \cite{CoSuTu}.
The method in \cite{CoSuTu}
based on the contraction mapping principle
and the Schauder fixed point theorem
goes through under the present assumptions on $a$ and $b$.
This gives the existence of solutions $z$, $w$ with the required
properties (i)-(iii) in the Sobolev class
$L^{1,p}(\D)$ for some $p > 2$.

The obtained solution is in $C^{1,\alpha}(\D)$ by elliptic
bootstrapping.
The $C^{1,\alpha}$ regularity up to the boundary follows
by the reflection principle \cite{IvSu} about the
totally real torus $b\D^2$.

The method of the proof of Theorem \ref{Riemann}
in \cite{CoSuTu} based on the Schauder
principle does not guarantee the continuous dependence
of solutions on the boundary conditions and coefficients of
the system (\ref{mainsystem}).
Instead, we use general results \cite{Mo}
on quasi-linear elliptic equations in the plane.
We reduce the system (\ref{mainsystem})
to the following single equation by eliminating
the variable $\zeta$:
\begin{eqnarray}
\label{Beltrami1}
w_{\overline z} +  a(z,w)  w_z  = b(z,w).
\end{eqnarray}
Indeed, using the first equation in (\ref{mainsystem}) we obtain
$$
w_{\bar\zeta}=w_z z_{\bar\zeta}+w_{\bar z}
{\bar z}_{\bar\zeta}=(a(z,w) w_z+w_{\bar z})
{\bar z}_{\bar\zeta}.
$$
Substituting the latter in the second equation of (\ref{mainsystem})
we obtain the equation (\ref{Beltrami1}).
This equation is equivalent to the original system (\ref{mainsystem})
because once the solution of (\ref{Beltrami1}) is found,
one can find $\zeta(z)$ from the linear Beltrami equation
$$
\zeta_{\bar z}+a(z,w(z))\zeta_z=0.
$$
The latter is in fact the first equation in (\ref{mainsystem})
written for the inverse function $\zeta(z)$.

The results of \cite{Mo} concern a
more general quasi-linear equation
\begin{eqnarray}
\label{CRgraph1}
w_{\overline z} + a_1(z,w) w_z + a_2(z,w)
\overline{w}_{\overline z}= b(z,w)
\end{eqnarray}
with a linear boundary condition
\begin{eqnarray}
\label{boundlinear}
\Re [\bar {G(\zeta)} w(\zeta)] = g(\zeta),\quad
\zeta \in b\D.
\end{eqnarray}
Let $\ind G$ denote the winding number of $G$.
It is called the index of the Riemann--Hilbert problem
(\ref{CRgraph1}--\ref{boundlinear}).
We assume that
\begin{itemize}
\item[(i)]
the coefficients $a_j$ and $b$ are $L^\infty$ in
$\D\times\cc$ and Lipschitz in $w$ uniformly in $z$;
\item[(ii)]
the ellipticity condition
$|a_1|+|a_2|\le a_0 < 1$
holds, here $a_0$ is constant;
\item[(iii)]  the functions $G$ and $g$ in (\ref{boundlinear})
are $C^\beta$ on the unit circle (for some $\beta > 0$ )
and $G\ne 0$;
\item[(iv)] $\ind G\ge0$.
\end{itemize}

Then the following result \cite{Mo} (pp. 335--351) holds.

\begin{e-pro}
\label{monakhov} Under the above assumptions (i)-(iv),
the boundary value problem
(\ref{CRgraph1}--\ref{boundlinear}) admits a solution $w$
in the Sobolev class $L^{1,p}(\D)$ for some $ p > 2$.
The solution is unique if it satisfies the conditions
$$
w(p_j) = 0, \quad
\Im w(1) = 0
$$
for some fixed points $p_j \in \D$, $j=1,\dots,\ind G$.
Furthermore,
the solution continuously depends on perturbations
of $G$, $g$ and $a_j$, $b$ in the $C^\beta$ and $L^\infty$
norms respectively.
\end{e-pro}

We point out that in \cite{Mo} this result is obtained
under substantially weaker regularity assumptions.

Return now to the equation (\ref{Beltrami1}) with the non-linear
boundary condition
\begin{eqnarray}
\label{boundnonlinear}
|w (z)|= r\quad
{\rm for}\quad
|z|= 1.
\end{eqnarray}
Set $w = re^{it}z^n e^u$.
Then the new unknown $u$ satisfies the equation
\begin{eqnarray}
\label{Beltrami2}
u_{\bar z} +
a(z,re^{it}z^n e^u) u_z =r^{-1}e^{-it}z^{-n}e^{-u}b(z,re^{it}z^ne^u)
- nz^{-1} a(z,re^{it}z^ne^u)
\end{eqnarray}
with the linear boundary condition
\begin{eqnarray*}
\Re u(z) = 0\quad
{\rm for}\quad
|z|= 1,
\end{eqnarray*}
and the index of the problem is equal to zero.
The coefficients of (\ref{Beltrami2})
are bounded and still Lipschitz in $u$ uniformly in $z$, $r$ and $t$.
By Proposition \ref{monakhov} the solution $u$ with $u(1)=0$
depends continuously on the parameters $r$, $t$
and the coefficients $a$ and $b$.
Hence the solution $(z,w)$ with $z(1)=1$ and $w(1)=re^{it}$
in Theorem \ref{Riemann} continuously depends on $r$, $t$
and the coefficients $a$ and $b$.
Since for $r=0$ we have $(z,w)=(\zeta,0)$,
then $(z_n,w_n)$ is homotopic to $(\zeta,0)$.
Finally, the conclusion (iv) follows by the uniqueness
and the fact that we can replace the condition
$w(1)=re^{it}$ by $w(\zeta_0)=re^{it}$ for a fixed
$\zeta_0\in b\D$.
Theorem \ref{Riemann} is now proved.

\section{Gluing $J$-holomorphic discs to real tori}

As the first application of Theorems \ref{removal} and \ref{Riemann}
we obtain a result on gluing pseudo-holomorphic discs to
real tori. In \cite{CoSuTu} we construct $J$-holomorphic discs
{\it approximately} attached to real tori.
We recall that for the usual complex structure,
Forstneri\v c and Globevnik \cite{FoGl}
gave constructions for approximately attaching
holomorphic discs to certain tori in $\cc^n$.
We improve one of the results of  \cite{CoSuTu}
here by constructing discs attached to tori exactly.

Throughout this section $(M,J)$ denotes
an almost complex manifold of complex dimension 2.
Let $f_0: \bar\D \to M$ be an immersed
$J$-holomorphic disc smooth on $\bar\D$. Fix $R > 1$ and
consider a family of $J$-holomorphic immersions
$h_z: R\bar\D \to M$ smooth on $R\overline{\D}$
and smoothly depending on the parameter $z \in \bar\D$.
Let $\Sigma$ denote the set of all critical points of the map
$H: \overline\D\times  R\overline\D\ni(z,w) \mapsto h_z(w)\in M$.
Suppose that  the following conditions hold.
\begin{itemize}
\item[(i)] For every $z \in \bar\D$ we have $h_z(0) = f_0(z)$ ;
\item[(ii)] For every $z \in \bar\D$ the disc $h_z$
is transverse to $f_0$ .
\item[(iii)] The map $H|_{(\D \times R\D)\setminus\Sigma}$
preserves the canonical orientations defined by $J_{st}$ and $J$ on
$\D \times R\D$ and $M$ respectively.
\end{itemize}
Define $\Lambda = H(b\D \times b\D)$. We call $\Lambda$
a real torus though it is even not immersed in general.
\begin{e-theo}
\label{tori}
Let $\Omega$ be a bounded strictly pseudoconvex
domain in $(M,J)$.
Let $f_0:\D \to \Omega$ be a $J$-holomorphic
immersion smooth in $\overline{\D}$.
Fix $R > 1$ and consider a smooth family
$h_z: R\overline\D \to \overline\Omega$,
$z \in \overline\D$   of
$J$-holomorphic immersions satisfying
(i), (ii), and (iii). Let $\Lambda$ be a real torus defined above.
Suppose that there exists  $c_0: 0 < c_0 < 1$ such that
$h_z(b\D) \subset b\Omega$ for every
$z: c_0 < \vert z \vert \leq 1$; in particular,
$\Lambda\subset b\Omega$.
Then there exists a continuous one-parameter family
of $J$-holomorphic discs $f^t:\bar\D\to\bar\Omega$
of some class $C^{1,\alpha}(\bar\D)$
such that $f^t(b\D) \subset \Lambda$,
$f^t(0) = f_0(0)$, and
$f^t$ is tangent to $f_0$ at the center.
The boundaries of the discs $f^t$ fill the whole torus
$\Lambda$.
\end{e-theo}

As we pointed out earlier, the condition (iii)
in Theorem \ref{removal} can be replaced by
(iii'): the set $\D^2 \backslash \Sigma$ is connected.

\proof The conditions (i)-(iii) allow to apply
Theorem \ref{removal} on $\D \times R\D$.
By Theorem \ref{removal} the pull-back $\tilde J=H^*(J)$
is well defined in the bidisc $\overline{\D}^2$ and has
a complex matrix of the form (\ref{Aab}), where
the coefficients $a(z,w)$ and $b(z,w)$ are
$C^\alpha$ in $z$ and $C^{0,1}$ in $w$ for some $0<\alpha<1$.
A map $\zeta \mapsto (z(\zeta),w(\zeta))$, $\zeta \in \D$,
is $\tilde J$-holomorphic if and only if it satisfies
the system (\ref{mainsystem}).

We extend the functions $a$ and $b$ to $\overline\D \times \cc$
keeping them $C^\alpha$ in $z$ and $C^{0,1}$ in $w$.
In the notation of Theorem \ref{Riemann}, fix an integer $n$
(depending on $c_0$) big enough and also fix $0\le t<2\pi$.
Let $(z_r, w_r)=(z_r^t, w_r^t)$, $0<r\le1$, be a family of solutions
of (\ref{mainsystem}) constructed by Theorem \ref{Riemann}
for these  $n$ and  $t$. We claim that the discs $(z_r, w_r)$
stay in $\overline{\D}^2$ for all $0<r\le1$. Then the disc
$f^t(\zeta)=H(z_1^t(\zeta),w_1^t(\zeta))$ satisfies
the conclusion of the theorem.

Arguing by contradiction, assume that the disc $(z_1, w_1)$
is not contained in $\overline{\D}^2$, that is,
$|w_1(\zeta)|>1$ for some $\zeta\in\D$.
Fix a constant $c_1$ such that  $c_0 < c_1 <1$.
By (ii) and (iii) of Theorem \ref{Riemann},
$|w_r(\zeta)|<1$ for $|z_r(\zeta)|< c_1$
if the above $n$ is chosen large enough.
If $r$ is small, then $|w_r|$ is also small.
Hence, there are $r$ and $\zeta\in\D$ such that
$c_0<|z_r(\zeta)|<1$ and $|w_r(\zeta)|=1$.
Then the disc $H(z_r,w_r)$ touches the strictly
pseudoconvex hypersurface $b\Omega$ from inside,
which is impossible.

The boundaries of the discs $f^t$ fill the
whole torus $\Lambda$ by (iv) in Theorem \ref{Riemann}.
Theorem \ref{tori} is proved.
\medskip

A special but important for applications situation arises
if the above map $H$ is an immersion i.e. the set $\Sigma$ is empty.
In this case it suffices to require that $H$ is defined just
on $b\D \times R\overline\D$. For convenience we state
the corresponding assertion explicitly.

\begin{e-cor}
Let $\Omega$ be a bounded strictly pseudoconvex
domain in $(M,J)$. Suppose that for some $R > 1$
a map $H:b\D \times R\overline \D \to \overline\Omega$ is
a smooth immersion satisfying the following assumptions:
\begin{itemize}
\item[(i')]
for every $z \in b\D$ the map $H(z,\bullet):R\D \to M$
is $J$-holomorphic;
\item[(ii')] the map $H(\bullet,0):\D \to M$ is
a $J$-holomorphic immersion.
\end{itemize}
Suppose that the torus $\Lambda = H(b\D \times b\D)$
is contained in $b\Omega$.
Then there exists a continuous one-parameter family
of $J$-holomorphic discs $f^t:\bar\D\to\bar\Omega$
of  class $C^{\infty}(\bar\D)$
such that $f^t(b\D) \subset \Lambda$,
$f^t(0) = f_0(0)$, and $f^t$ is tangent to $f_0$ at the center.
The boundaries of the discs $f^t$ fill the whole torus
$\Lambda$.
\end{e-cor}

\proof It is shown in \cite{CoSuTu} that after a suitable
reparametrization $\zeta\mapsto e^{i\sigma(z)}\zeta$
of the discs $h_z$, the map $H$ defined above extends
smoothly to the whole bidisc $\overline{\D}^2$ as an immersion such
that  the above assumption (i') is satisfied. Moreover,
we choose the discs $h_z$ in such a way that
$h_z(b\D)\subset b\Omega$ for $c_0<|z|\le1$, where
$c_0$ is close to 1. Now we can apply Theorem \ref{tori}.
This completes the proof.

\bigskip

{\bf Remarks.}

1. We use in the proof that the defining
function $\rho$ is strictly plurisubharmonic in a neighborhood
of $b\Omega$, not necessarily on all of $\Omega$.

2. A version of Theorem \ref{Riemann} still holds with the same
proof if $w$ is vector valued and satisfies the equations
$(w_j)_{\overline\zeta}= b_j(z,w)\overline z_{\overline\zeta}$
and boundary conditions
$|w_j(\zeta)|=r$, $|\zeta|= 1$.
However, in order to prove the corresponding version of
Theorem \ref{tori} in higher dimension
we need $J$-complex hypersurfaces in $M$.
They generally do not exist in higher dimension unless
$J$ is integrable.
\medskip

In Theorem \ref{tori} we suppose that $h_z(b\D) \subset b\Omega$
for $c_0 < \vert z \vert \leq 1$. In some applications
the stronger condition $H(\overline\D \times b\D) \subset b\Omega$
holds i.e. $h_z(b\D) \subset b\Omega$ for all $z \in \overline\D$.
In this special case one
can construct discs with additional properties.

\begin{e-theo}
\label{toridouble}
Let $\Omega$ be a bounded strictly pseudoconvex
domain in $(M,J)$. Let a smooth map
$H: \bar\D^2 \to \bar\Omega$
satisfies the assumptions of Theorem \ref{removal}.
Let the map $z\mapsto f_0(z):= H(z,0)$
be $J$-holomorphic and smooth in $\bar\D$.
Suppose that $H$ is an immersion near $\overline \D \times b\D$
and $H(\overline\D \times b\D)\subset b\Omega$.
Set $\Lambda := H(b\D \times b\D)$.
Fix finitely many points $z_j \in \D$, $j=1,\dots,s$
and integers $m_j \ge 0$.
Finally suppose the points $f_0(z_j)$ are distinct
non-critical values of $H$.
Then there exists a continuous one-parameter family
of $J$-holomorphic discs $f^t:\bar\D\to\bar\Omega$
of some class $C^{1,\alpha}(\bar\D)$
such that $f^t(b\D)\subset\Lambda$ and
$f^t$ is tangent to $f_0$  at the points $f_0(z_j)$
with orders $m_j$ respectively.
The boundaries of the discs $f^t$ fill the whole torus
$\Lambda$.
The family $f^t$ depends continuously on perturbations
of $H$ and $J$.
 \end{e-theo}

The proof is similar to the proof of Theorem \ref{tori}.
The only difference is that we apply Proposition \ref{monakhov}
instead of Theorem \ref{Riemann}. We also point out that
here we do not need to assume that $H$ is defined on
$\overline \D \times R\overline\D$ for some $R > 1$.
Indeed, the assumption on the map $H$ to be an immersion near
$\overline \D \times b\D$ implies that  the intersection of
$\Sigma$ with every set $\{ z \} \times \D$,
$z \in \overline\D$ is compact and Theorem \ref{removal}
can be applied.
\medskip

{\bf Remarks.}

1. If we assume in Theorems \ref{tori} and \ref{toridouble}
that $H$ is an immersion on $b\D^2$, then we can prove
them without using the continuous dependence statement
in Theorem \ref{Riemann}. Moreover, the resulting disc will
be $C^\infty$ smooth up to $b\D$ by the reflection principle
\cite{IvSu} about the totally real torus $\Lambda$.

2. Theorems \ref{tori} and \ref{toridouble}
are new also in the case where the manifold $M$ is
the Euclidean space $\cc^2$ with the standard complex
structure.

\section{Appendix. Proof of Theorem \ref{removal}
in the integrable case}

Assume that in the setting of Theorem \ref{removal}
the structure $J'$ is integrable.
Then without loss of generality $M=\cc^2$ and $J'=J_{st}$.
By Lemma \ref{localization} we can restrict to
a local problem near the origin.
Let $H = (h_1,h_2)$.
Here the functions $h_j$ are $C^k$ in $(z,w)$ and holomorphic
in $w$; $k\ge1$ is not necessarily an integer.

Since every map $w \mapsto H(z,w)$ is an immersion, then
$(h_2)_w(0,0) \neq 0$.
Then we perform a $C^k$ diffeomorphic change of variables
$(z,w)\mapsto(z,h_2(z,w))$ in a neighborhood of $(0,0)$
keeping the same notation for the new variables.
Note that the change of variables preserves the
coordinate lines $z=c$, hence it reduces the problem
to the case in which
$$
H(z,w) = (h(z,w),w).
$$
Then the singular set has the form
$\Sigma=\{(z,w): |h_z|=|h_{\bar z}|\}$.
Since $H$  is orientation preserving, then
$|h_z|\ge|h_{\bar z}|$.
Since $h$ is holomorphic in $w$,
then so are $h_z$ and $h_{\bar z}$.
By the condition (ii), for every fixed $z$, the function
$h_z(z,\bullet)$ is not identically equal to zero.
Let $a=-h_{\bar z}/h_z$.
Then for every fixed $z$ the function $a(z,w)$
is holomorphic in $w$ with isolated singularities
at zeros of $h_z$, and $|a|\le 1$.
Then the singularities are removable, and by the maximum
principle and condition (ii) again we have $|a|<1$.
Then we immediately obtain
$\Sigma=\{(z,w): h_z(z,w)=0\}$, hence
$\Sigma$ is discrete for fixed $z$.

Let $J=H^*J_{st}$.
By the transformation rule (\ref{PDEfg4})
the complex matrix $A$ of $J$
in the complement of $\Sigma$ has the form
\begin{eqnarray*}
A =\left(
\begin{array}{cll}
a &  0\\
0&  0
\end{array}
\right),
\end{eqnarray*}
where $a$ is defined above.
Since $|a|< 1$, then $A$, whence $J$ extends
to $\Sigma$.
The structure $J$ is $C^{k-1}$-smooth in $(z,w)$
by the Cauchy integral formula in $w$:
$$
a(z,w) =  \frac{1}{2\pi i} \int_\gamma
\frac{a(z,\tau)\,d\tau }{\tau -w}
$$
where $\gamma$ is a simple path in the complex
$w$-plane so that $h_z(z,w)\ne0$ for all $w\in\gamma$ and $z$
in a small open set.
The proof is complete.

\end{document}